\documentclass{elsart}
\usepackage{graphicx}
\begin{document}
\begin{frontmatter}
\title{Hilbert series of algebras associated to directed graphs and order homology}

\author{Vladimir Retakh}
\address{Department of Mathematics, Rutgers University,
Piscataway,
NJ 08854-8019, USA}
\ead{vretakh@math.rutgers.edu}

\author{Shirlei Serconek}
\address{IME-UFG,
CX Postal 131,
Goiania - GO,
CEP 74001-970,  Brazil}
\ead{serconek@math.rutgers.edu}

\author{Robert Wilson}
\address{Department of Mathematics, Rutgers University,
Piscataway,
NJ 08854-8019, USA}
\ead{ rwilson@math.rutgers.edu}

\begin{abstract}
We give a homological interpretation of the coefficients of the Hilbert series
for an algebra associated with a directed graph and its dual algebra. This allows
us to obtain necessary conditions for Koszulity of such algebras in terms of homological
properties of the graphs. We use our results to construct algebras with a prescribed
Hilbert series.
\end{abstract}


\end{frontmatter}
\maketitle

\newcommand{\F}{\frak F}
\linespread{1.6}  

\addtolength{\textwidth}{1in}
\addtolength{\hoffset}{-.5in}

\def\pbmatrix#1{\begin{bmatrix} #1 \end{bmatrix}}
\def\reals{ {\mathbb R}}
\def\complex{ {\mathbb C}}
\def\natural{ {\mathbb N}}
\def\rats{ {\mathbb Q}}
\def\integer{ {\mathbb Z}}

\def\ff{ {\Bbb F}}
\def\ba{ {\bf a}}
\def\bb{ {\bf b}}
\def\bc{ {\bf c}}
\def\bu{ {\bf u}}
\def\bv{ {\bf v}}
\def\bw{ {\bf w}}

\parindent=0cm

\newtheorem{example}{Example}[subsection]
\newtheorem{lemma}[example]{Lemma}
\newtheorem{proposition}[example]{Proposition}
\newtheorem{theorem}[example]{Theorem}
\newtheorem{corollary}[example]{Corollary}


\section*{Introduction}
\label{sec: 0}

In \cite {GRSW,GGRW1,RSW1,RSW2,RSW3,RSW4,RW2} we introduced and studied certain associative
noncommutative algebras $A(\Gamma )$ defined by layered graphs
$\Gamma$ (or ranked posets) and their generalizations. The algebras $A(\Gamma)$ are
related to factorizations of polynomials with noncommutative coefficients and we called
them {\it splitting algebras}. An important example of such algebras is the algebra $Q_n$
defined by the Boolean lattice of subsets of a finite set (see {\cite {GRW}}) and
related to the theory of noncommutative symmetric functions {\cite{GKLLRT}}. For
{\it homogeneous} layered graphs \cite {RSW1} the algebras $A(\Gamma)$ are quadratic and one can
construct their quadratic dual algebras $A(\Gamma)^!$.

It turns out that algebraic properties of $A(\Gamma)$ and $A(\Gamma)^!$
are closely related to homological properties of $\Gamma$. When a layered
graph $\Gamma$ is defined by a ``good" regular
cell complex the algebras $A(\Gamma)$ and $A(\Gamma)^!$ are Koszul if
and only if all intermediate (order) cohomologies of $\Gamma$ are trivial ( see {\cite {CPS,RSW4}}).

Recall that if a quadratic algebra $A$ is Koszul then it is {\it numerically} Koszul,
i.e. the Hilbert series (or the graded dimension) $h(A,\tau)$ of $A$ and
the Hilbert series $h(A^!,\tau)$ of its dual algebra $A^!$ satisfy
the identity
\begin{equation} \label{eq:numkoz}
h(A,\tau)h(A^!, -\tau)=1.
\end{equation}
In this paper, motivated by the results from \cite {CPS,RSW4} and identity (\ref{eq:numkoz})
we present a homological interpretation of coefficients of the Hilbert series for
the algebras $A(\Gamma)$ and $A(\Gamma)^!$. In fact, instead of $A(\Gamma)^!$
we are working with a simpler algebra, the algebra $B(\Gamma), $
the graded algebra associated with a natural filtration on $A(\Gamma)^!$.
The algebras $A(\Gamma)$, $A(\Gamma)^!$, $B(\Gamma)$ are Koszul if and only if at least
one of these algebras is Koszul. Also, $h(A(\Gamma)^!,\tau)=h(B(\Gamma),\tau)$.
Note that the algebras $B(\Gamma)$ are defined for any directed graph without any extra assumptions.

Our homological interpretation of the coefficients of the Hilbert series for the algebras
$A(\Gamma)$ and $B(\Gamma)$ allows us to obtain conditions for their
numerical Koszulity and also to construct
algebras with prescribed Hilbert series defined, for example, by  palindromic
polynomials.
Our formulas look particularly simple when certain posets associated with $\Gamma$
are Cohen-Macaulay.

There are several ways to associate an algebra to a directed graph or a poset, and  there are
known connections between properties of such algebras and topological structures of graphs and posets
(see \cite {Bjo}). The most famous example is the incidence algebra of a finite poset
described in \cite {St}, Section 3.6. There is a certain resemblance between
our algebra $A(\Gamma)^!$ and the incidence algebra defined by the graph $\Gamma$ but
our results are quite different.

The paper is organized in the following way. In Section 1 we recall basic facts
about graphs, posets and their (co)homologies. Section 2 contains the definition
of the algebras $A(\Gamma)$, $A(\Gamma)^!$ and $B(\Gamma)$. Our main theorem on homological
description of the coefficients in the Hilbert polynomials for the algebras $B(\Gamma)$ and its
corollaries are formulated in Section 3. Section 4 is devoted to numerical Koszulity
of the algebras $A(\Gamma)$ and $B(\Gamma)$. Section 5 contains a number of examples including examples
of algebras with Hilbert series equal to $P(-\tau)^{-1}$ where
$P(\tau)$ is a palindromic polynomial.
Calabi-Yau algebras also have Hilbert series defined by palindromic
polynomials (see {\cite {G}}).

During preparation of this paper Vladimir Retakh and Robert Wilson were supported by an
NSA grant.

\section{Partially ordered sets and directed graphs}

\subsection{Partially ordered sets and their homology}

Let $P$ be a partially ordered set (poset).
Define an $i$-chain, i.e. a chain of length $i$, in $P$ to be an $(i+1)$-tuple
$${\mathbf x} = (x_0,...,x_i)$$
of elements of $P$ with
$$x_0 < x_1 < ... < x_i.$$
Let $Ch_i(P)$ denote the set of $i$-chains.
The M\"obius function, $\mu$ is defined by
$$\mu(x,y) = c_0-c_1+c_2-\dots$$
where $c_i$ is the
number of $i$-chains ${\mathbf x} = (x_0,...,x_i)$ with
$x = x_0, y = x_i$.

We say that $u\in P$ covers $v\in P$ (and write $v\prec u$) if $v<u$ and there is no element between
$u$ and $v$. If $P$ is finite, denote by $l(P)$ the maximal length of a chain
in $P$.

Denote by $\hat P$ the poset obtained from $P$ by adding the minimal element $\hat 0$
and the maximal element $\hat 1$. We write $\mu (P)$ for $\mu (\hat 0, \hat 1)$.

Suppose that $P$ is a lattice, i.e. for any two elements $x,y\in P$
their least upper bound $x\vee y$ and the greatest lower bound
$x\wedge y$ are defined. A lattice $P$ of finite length is a lower semimodular lattice
if for any two elements $x,y\in P$,
if $x\vee y$ covers both $x$ and $y$ then $x$ and $y$ both cover $x\wedge y$.

Let $F$ be a field and $P$ be a finite poset. If $i \ge -1$, denote by $C_i(P;F)$
the free $F$-module on the set of
$i$-chains ${\mathbf x}$ of $P$. The empty set is a $(-1)$-chain, so we
identify $C_{-1}(P;F)$ with $F$. Set $C_n(P;F)=0$ if $n<-1$.

If ${\mathbf x} = (x_0,...,x_i)$ is an $i$-chain and $0 \le l \le i$, define $g^l({\mathbf x})$ to be the
$(i-1)$-chain
$$(x_0,...,x_{l-1},x_{l+1},...,x_i).$$
Define a map
$d_i: C_i(P;F)\rightarrow C_{i-1}(P;F)$ by linearly extending

\begin{equation} \label{eq:boundary}
d_i({\mathbf x})=\sum _{l=0}^i(-1)^l g^l({\mathbf x})
\end{equation}

when $i-1$-chains exist, and setting $d_i=0$ otherwise. It is easy to check
that $d_id_{i+1}=0$, and one can define {\it reduced order homology groups} of $P$
with coefficients in $F$:
$$\tilde H_i(P;F)=Ker\ d_i/Im\ d_{i+1}.$$

By construction, $\tilde H_i(P;F)=0$ for $i<-1$ and $i>l(P)$. Also, $\tilde H_{-1}(P;F)=0$
if and only if $P$ is nonempty, and $\tilde H_0(P;F)=0$ if and only if $P$ is connected.

Denote by $C^i(P;F)$ the vector space dual to $C_i(P;F)$.  For an $i$-chain $\mathbf x$ let
${\mathbf x}^*$ denote the element of $C^i(P;F)$ defined by
$${\mathbf x}^*({\mathbf y}) = \delta_{{\mathbf x}, {\mathbf y}}$$
for all $i$-chains $\mathbf y$.  Define $\partial^i$ to be
the linear map dual to $d_i$. One can then define  the cohomology groups
$$\tilde H^i(P;F) = Ker \ \partial ^{i+1}/Im \ \partial ^i.$$

The M\"obius function is the Euler characteristic of reduced order homology:
$$
\mu (P)=\sum _{i=-1}^{l(P)}(-1)^i\dim \ \tilde H_i(P;F).
$$

Corresponding to a poset $P$ there is a simplicial complex $\Delta(P)$.  The vertices of $\Delta(P)$ are the elements of $P$ and the $i$-faces of $\Delta(P)$ are the $i$-chains of $P$.

Recall, that a poset $P$ is {\it Cohen - Macaulay} if for each open interval
$(x,y)$ in $P$ all homology groups $\tilde H_i((x,y);F)$ are trivial for
$i < \dim \ \Delta((x,y))$. Any lower semimodular lattice is
Cohen-Macaulay.

From now on we will consider $F$  to be fixed and  write 
$\tilde H_i(P)$ for  $\tilde H_i(P;F)$ and
$\tilde H^i(P)$ for  $\tilde H^i(P;F)$.

\medskip

\subsection{Layered graphs}

Let $\Gamma =(V,E)$ be a directed graph (quiver) where $V$ is the set of vertices and $E$ is the set of edges.
For any $e\in E$ denote by $t(e)$ the tail of $e$ and by $h(e)$ the head of $e$.
A path $\pi$ in $\Gamma$
is a sequence of edges $\pi=(e_1,e_2,\dots, e_k)$ such that $h(e_i)=t(e_{i+1})$ for $i=1,2,\dots, k-1$. We call
$t(e_1)$ the tail of $\pi$ and denote it by $t(\pi)$ and we call
$h(e_k)$ the head of $\pi$ and denote it by $h(\pi)$.

Assume that $V=\coprod_{i=0}^N V_i$. We call $N$ the {\it height} of the graph
and we say that $i$ is the {\it level} of $v$ and write $|v|=i$ if $v\in V_i$.
We call $\Gamma $ a {\it layered} graph if
$|t(e)| = |h(e)| + 1$ for any edge $e \in E$. 

For any directed graph $\Gamma =(V,E)$ there is a corresponding partially ordered set.
The elements of the poset are the vertices $v\in V$. We say that $u>v$ if and only if
there is a directed path from $u$ to $v$. By abuse of notation, we denote the poset
by the same letter $\Gamma$. Therefore, we may talk about M\"obius functions and
(co)homologies
of directed graphs by considering them as posets.

Obviously, for a layered graph $u$ covers $v$ if and only if there is an edge
going from $u$ to $v$.

We say that vertices $v,v'$ of the same level $i>0$ are connected by a {\it down-up}
sequence if there exist vertices $v=v_0,v_1,v_2,\dots ,v_k=v'\in V_i$ and
$w_1,w_2,\dots, w_k\in V_{i-1}$ such that $w_j < v_{j-1}, v_j$ for $j=1,2,\dots , k$.
According to {\cite {RSW1}}, the layered graph $\Gamma $ is {\it uniform} if for any
pair of edges $e,e'\in E$ with a common tail, $t(e)=t(e')$, their heads $h(e), h(e')$
are connected by a down-up sequence $v_0,\dots, v_k$ such that $v_i < t(e)$,
$i=0,1,\dots, k$.

Let $\Gamma = (V,E)$ be a layered graph with $ V = \coprod_{i=0}^N V_i$ and $V_i = \emptyset$ unless $r \le i \le r+L-1$.
Note that this implies that $Ch_L(\Gamma) = \emptyset$ and that if
${\mathbf x} = (x_0,...,x_{L-1})$ is an $(L-1)$-chain then $x_i \in V_{r+i}$ for $0 \le i \le L-1$.

For $0 \leq l \leq L-1$ set
$$g^l(Ch_{L-1}(P)) =  Ch_{L-2}^{(l)}(P)$$
and
$$g^l(C_{L-1}(P)) =  C_{L-2}^{(l)}(P).$$

Let $C^{(l),L-2}(P)$  denote the set of all linear functions $f:C_{L-2}^{(l)}(P) \rightarrow F$. Now
$Ch_{L-2}(P)$ is the disjoint union of the $Ch_{L-2}^{(l)}(P)$  and so
 $$C^{L-2}(P) = \bigoplus _{l=0}^{L-1} C^{(l),L-2}(P)$$
where we extend $f \in C^{(l),L-2}(P)$ to a function on $C_{L-2}(P)$ by setting $f(C_{L-2}^{(m)}(P)) = 0$
for $m \ne l.$

\subsection{Examples of layered graphs}

Two families of layered graphs are particularly interesting.
\begin{example}
   Let $\Gamma = (V,E)$
where $V = \coprod_{i=0}^N V_i$ is a layered graph.  We say that $\Gamma$ is a {\it {complete layered graph}} if for every $i, 1 \le i \le N$,
and every $v \in V_i, w \in V_{i-1}$ there is an edge from $v$ to $w$.
Clearly a complete layered graph is determined up to isomorphism by the integers
$|V_N|,...,|V_1|, |V_0|$.    Let $m_N,...,m_0$ be positive integers.
Denote by ${\mathbf C}[m_N,...,m_0]$ the complete layered graph with
$|V_i| = m_i$ for $0 \le i \le N$.
\end{example}

\begin{proposition}
$$\dim H^{i}({\mathbf C}[m_N,...,m_0]) = 0$$
if $-1 \le i \le N-1$ and
$$\dim H^{N}({\mathbf C}[m_N,...,m_0]) = (m_N - 1)(m_{N-1} -1)...(m_0 - 1).$$
\end{proposition}

\noindent {\bf Proof:}
The result is clear if $N = 0$ or $-1$.  Hence assume $N > 1$. We first show that the partially ordered set $V$ is lexicographically
shellable (in the sense of Definition 2.2 of \cite{Bjo}).
For each $i, 0 \le i \le N$ choose $x_i' \in V_i$.  If $e \in E_i$,
define $\lambda(e) = 1$ if $t(e) = x_i'$ and $\lambda(e) = N+1-i$ otherwise.
Then if $[x,y]$ is an interval of $V$ with $x \in V_i, y \in V_j$ we
see that $x < x_{i+1}'< ... <x_{j-1}' < y$ is the unique rising
unrefinable chain from $x$ to $y$, so this labeling is an $R$-labeling.
Furthermore if $x \prec z \le y$ and $z \ne x_{i+1}'$ we
have $\lambda(x,x_{i+1}') = 1 < N+1-i = \lambda(x,z).$
Thus this labeling is an $L$-labeling and so $V$ is lexicographically
shellable. Then  by Theorem 3.2 of \cite {Bjo}  $V$ is Cohen-Macaulay, giving the first statement.

Now
$$\dim C^i({\mathbf C}[m_N,...,m_0]) = \dim C_i({\mathbf C}[m_N,...,m_0])=$$
$$|Ch_i({\mathbf C}[m_N,...,m_0])|=\sum _{|J|=i+1, J\subseteq \{0,\dots ,N\}}\prod _{j\in J}m_j.$$

Thus, by the first part of the proposition and the Euler-Poincare principle,
$$\dim H^N({\mathbf C}[m_N,...,m_0]) =
(-1)^N\sum _{i=-1}^N (-1)^i \dim H^i({\mathbf C}[m_N,...,m_0]) = $$
$$(-1)^N\sum _{i=-1}^N(-1)^i\dim C^i({\mathbf C}[m_N,...,m_0])) =
 (-1)^N\sum _{i=-1}^N(-1)^i\sum _{|J|=i+1}\prod _{j\in J}m_j = $$
$$ = \sum _{J\subseteq \{0,\dots ,N\}}(-1)^{N + 1 -|J|}\prod _{j\in J}m_j=\prod _{i=0}^N(m_i-1).$$

\vskip 20 pt
\begin{example}   Let $\Theta_N$ denote the Hasse graph of the (Boolean) lattice of
all subsets of $\{1,...,N\}$.  Thus the vertices of $\Theta_N$ are the subsets of $\{1,...,N\}$, the level of
$Y \subseteq \{1,...,N\}$ is $|Y|$,  and there is an edge for $Y \subseteq \{1,...,N\}$ to $Z \subseteq \{1,...,N\}$ if and only if $Y \supseteq Z$ and $|Y| = |Z| + 1$.
\end{example}
If $\Gamma$ is any layered graph, $a \in V(\Gamma)$, and $i \le |a|$, define $\Gamma_{a,i}$ to be the subgraph induced by the set of vertices
$$\{w \in V(\Gamma)|a > w \ \ \textrm{and} \ \ |a| - |w| \le i-1\}.$$
If $\Gamma$ has unique maximal vertex $v$ denote $\Gamma_{v,i}$ by $\Gamma_i$.  Clearly
$(\Theta_N)_{a,i}$ is isomorphic to $\Theta_{|a|,i}.$

We will need the following results about the order homology of $\Theta_N$ later.

\begin{proposition}  If $2 \le i \le N$ and $1 \le j \le i-3$ then
 $$H^j(\Theta_{N,i}) = (0).$$
\end{proposition}
 \noindent {\bf Proof:} Note that any edge $e$ in $\Theta_{N,i}$ has tail $S$ and head $S \setminus \{a\}$
for some $S \subseteq \{1,...,N\}$ and some $a \in S$.  Define $\lambda(e) = a.$  Then $\lambda$ is
an $L$-labeling in the sense of Definition 2.2 of {\cite {Bjo}} and so,
$\Theta_{N,i}$ is lexicographically shellable and hence (cf. \cite {Bjo}) $\Theta_{N,i}$ is Cohen-Macaulay, giving the result.

\begin{proposition} If $2 \le i \le N $ we have
$$\dim \ H^{i-2}(\Theta _{N,i}) ={{N-1}  \choose {i-1}}.$$
\end{proposition}

\noindent {\bf Proof:} If $A \subseteq \{1,...,N\}$, let $A^c$ denote $\{1,...,N\}\setminus A.$
Let $N[i-1]$ denote the set of all ordered $(i-1)$-tuples ${\mathbf x} = (x_0,...,x_{i-2})$ of distinct
elements of $\{1,...,N\}$.  Note that $Sym_{i-1}$, the symmetric group on $\{0,...,i-2\}$, acts on
$N[i-1]$ by permuting subscripts.  For $1 \le l \le i-2$ let $\sigma_l$ denote the transposition interchanging
$l-1$ and $l$.

If ${\mathbf x} = (x_0,...,x_{i-2}) \in N[i-1]$ and $0 \le m \le i-2$ define
$$b_m({\mathbf x}) = \{1,...,N\} \setminus \{x_m,...,x_{i-2}\}$$
and
$$b({\mathbf x}) = (b_0({\mathbf x}),...,b_{i-2}({\mathbf x})).$$
Then
$$b:N[i-1] \rightarrow Ch_{i-2}(\Theta_{N,i})$$
is a bijection and
$\{b({\mathbf x})^*| {\mathbf x} \in N[i-1] \}$ is a basis for $C^{i-2}(\Theta_{N,i}).$

Since $b$ and $g^{l}$ are both surjective, the composition
$$g^{l}b: N[i-1] \rightarrow Ch_{i-3}^{(l)}(\Theta_{N,i})$$
is surjective.  Furthermore, if $1 \le l \le i-2$ and ${\mathbf x}, {\mathbf y} \in N[i-1]$
we see that
$$ g^{l}b({\mathbf x}) = g^{l}b({\mathbf y})$$
if and only if
$${\mathbf y} \in \{{\mathbf x}, \sigma_l{\mathbf x}\}.$$
Also, if ${\mathbf x} = (x_0,...,x_{i-2}), {\mathbf y} = (y_0,...,y_{i-2}) \in N[i-1]$, then
$$g^{0}b({\mathbf x}) = g^{0}b({\mathbf y})$$
if and only if
$$(x_1,...,x_{i-2}) = (y_1,...,y_{i-2}).$$

Now $$d_{i-2}(b({\mathbf y})) = \sum_{l=0}^{i-2} (-1)^lg^{l}b({\mathbf y)}$$
and so
$$\partial^{i-1}((g^{l}b({\mathbf y}))^*) = (-1)^l\sum_{g^lb({\mathbf x}) = g^{l}b({\mathbf y})}  b({\mathbf x})^*.$$
Hence, if $1 \le l \le i-2$,
$$\partial^{i-2}(C^{(l),i-3}(\Theta_{N,i})) = span \{b({\mathbf x})^* + b(\sigma_l{\mathbf x})^*|{\mathbf x} \in N(i-1)\} \subseteq C^{i-2}(\Theta_{N,i}).$$
Thus $\sigma_l$ acts as the identity on $\partial^{i-2}(C^{(l),i-3}(\Theta_{N,i}))$  so setting
$$\Xi =\sum_{\sigma \in Sym_{i-1}} (sgn  \ \sigma)\sigma$$ we see that
$$\Xi \partial^{i-2}(C^{(l),i-3}(\Theta_{N,i}))=(0)$$
for $1 \le l \le i-2.$
Furthemore,
$$\partial^{i-2}(C^{(0),i-3}(\Theta_{N,i-1})) = $$
$$span \ \{\sum_{y_0 \in \{x_1,...,x_{i-2}\}^c} b(y_0,x_1,...,x_{i-2})^*|{\mathbf x} = (x_0,...,x_{i-2}) \in N[i-1]\}.$$

Consequently, if $x_0 = 1$, then
$$b({\mathbf x})^* + \sum_{1 < y_0 \in \{x_1,...,x_{i-2}\}^c} b(y_0,x_1,...,x_{i-2})^* \in \partial^{i-2}(C^{i-3}(\Theta_{N,i})).$$
Thus if
$$D = span \  \{b({\mathbf x})^*|{\mathbf x} \in N[i-1], 1 < x_0 < x_1 < ... < x_{i-2} \}$$
we have
$$C^{i-2}(\Theta_{N,i}) = \partial^{i-2}(C^{i-3}(\Theta_{N,i})) + D.  $$


Define a linear map
$$\Psi: C^{i-2}(\Theta_{N,i}) \rightarrow C^{i-2}(\Theta_{N,i})$$
by
$$\Psi(b({\mathbf x})^*) = b({\mathbf x})^* \ \ \textrm{if} \ \ 1 \notin \{x_0,...,x_{i-2}\}$$
and
$$\Psi(b({\mathbf x})^*) = - \sum_{y_l \in \{x_0,...,x_{i-1}\}^c} b(x_0,...,x_{l-1},y_l,x_{l+1},...,x_{i-2})^*  \ \ \textrm{if}  \ \ x_l = 1.$$

Then $\Psi\sigma = \sigma\Psi$ for all $\sigma \in Sym_{i-1}$.  Thus
$$\Xi\Psi (\partial^{i-2}C^{(l),i-3}(\Theta_{N,i})) = (0)$$
for $1 \le l \le i-2.$

Now if ${\mathbf x } \in N[i-1]$ and $ 1 \notin \{x_0,...,x_{i-2}\}$ then $\Psi \partial^{i-2}(g^{0}b({\mathbf x})^*) =
\Psi(b(1,x_1,...,x_{i-2})^*) + \sum_{1 < y_0  \in \{x_1,...,x_{i-2}^c\}}
 b(y_0,x_1,...,x_{i-2})^* = 0.$
Also if
$1 = x_l$ with $1 \le l \le i-2$ we have
$$\Psi\partial^{i-2}((g^{0}b(\mathbf x))^*) = -\sum_{y \ne z, y,z \in \{x_0,...,x_{i-2}\}^c} b(y,x_0,...,x_{l-1},z,x_{l+1},...,x_{i-2})^*.$$
Since this is invariant under  the transposition that interchanges $0$ and $l$,
it is annihilated by $\Xi.$  Thus
$$\Xi\Psi\partial^{i-2}(C^{(0),i-3}(\Theta_{N,i})) = (0)$$
and so
$$\Xi\Psi\partial^{i-2}(C^{i-3}(\Theta_{N,i})) = (0).$$

Now $\Psi$ acts as the identity on $D$ and $\Xi$ is injective on $D$.  Thus we have
$$C^{i-2}(\Theta_{N,i}) = \partial^{i-2}(C^{i-3}(\Theta_{N,i})) \oplus D.  $$
Since
$$\dim  \ D = {{N-1} \choose {i-1}},$$
the proof of the proposition is complete.

\section{Algebras associated with layered graphs}
\subsection{The algebras $A(\Gamma)$}

Following {\cite {GRSW}} we construct now an algebra $A(\Gamma)$ associated with a layered graph $\Gamma$.
The algebra $A(\Gamma)$ is generated over the field $F$ by
generators $e\in E$ subject to the following relations. Let $t$ be a formal parameter commuting with edges $e\in E$.
Any two paths $\pi=(e_1,e_2,\dots,e_k)$ and
$\pi'=(f_1,f_2,\dots,f_k)$ with the same tail and head define the relation

\begin{equation} \label{eq;polynomial}
(t-e_1)(t-e_2)\dots (t-e_k)=(t-f_1)(t-f_2)\dots (t-f_k).
\end{equation}

In fact, relation (\ref{eq;polynomial}) is equivalent to $k$ relations
$$e_1+e_2+\dots +e_k=f_1+f_2+\dots +f_k,$$
$$\sum _{i<j}e_ie_j=\sum _{i<j}f_if_j, $$
$$\dots$$
$$e_1e_2\dots e_k=f_1f_2\dots f_k.$$

We call $A(\Gamma)$ the {\it splitting} algebra associated with graph $\Gamma$. The terminology is justified by
the following considerations.
Assume that there are only one vertex $*$ of the minimal level $0$ and only one vertex $x$ of the maximal level $N$,
and that for any edge $e$ there exists
a path $\theta =(e_1,e_2,\dots ,e_N)$ containing $e$ from $x$ to $*$. Set
$$P(t)=(t-e_1)(t-e_2)\dots (t-e_N).$$

Then $P(t)$ is a polynomial over $A(\Gamma)$
and any path from the maximal to minimal vertex corresponds to a factorization of $P(t)$ into a product of
linear factors.

\subsection{Dual algebras for uniform layered graphs}

Let $\Gamma =(V,E)$ be a layered graph. We assume that
the layered graph $\Gamma$ has exactly one minimal vertex so that
for any vertex $v\in V_i$, $i>0$  there is a path $\pi _v=(e_1,e_2,\dots, e_i)$ from $v$ to
the minimal vertex. In this case the splitting algebra $A(\Gamma)$ is defined by a set of homogeneous
relations of order $2$ and higher.

 It was proved in {\cite {RSW1}} that if the graph $\Gamma $ is uniform
then the splitting algebra $A(\Gamma)$ is quadratic, i.e. defined by relations of order $2$.

Recall that for a quadratic algebra $A$ over a field $F$ there is a notion of the dual quadratic algebra
$A^!$. To define $A^!$, denote by $W$ the $F$-span of the generators of $A$ and by $R\subset W\otimes W$
the linear space of relations of $A$. Denote by $W^*$ the dual space of $W$ and by
$R^{\perp}$ the annihilator of $R$ in $W^*\otimes W^*$. The algebra $A^!$ is the quadratic
algebra defined by generators $W^*$ and relations  $R^{\perp}$. It is well-known (see, for example,
{\cite {PP}}) that an algebra $A$ is Koszul if and only if its dual algebra $A^!$ is. In this case their
Hilbert series are connected by (\ref{eq:numkoz}).

\medskip \noindent
Assuming that a layered graph $\Gamma$ is uniform one can describe the dual algebra
$A(\Gamma )^!$ in terms of vertices and edges of the graph (see {\cite {RSW3}}). We describe
now a slightly different algebra $B(\Gamma )$.

There is a natural filtration on $A(\Gamma )$ defined by the ranking function $|\cdot |$.
The corresponding associated graded algebra is also quadratic.
Its dual algebra $B(\Gamma )$
can be described in the following way (see {\cite {CPS}}). Set $V_+=\coprod_{i>0}V_i$. For any $v\in V_+$ let $S(v)$ be the set of all vertices
$w\in V$ such that there is an edge going from $v$ to $w$.

\begin{theorem} The algebra $B(\Gamma )$ is generated by vertices $v\in V_+$ subject to the relations:

\noindent i) $\displaystyle u\cdot v=0$ if there is no edge going from $u$ to $v$;

\noindent ii) $\displaystyle v\cdot \sum _{ w\in S(v)}w=0$.
\end{theorem}

If the set of vertices $V$ is finite, the algebra $B(\Gamma )$ is finite-dimensional.
The algebras $B(\Gamma)$ were studied in {\cite {CPS}}.
According to the general theory $h(A(\Gamma)^!,\tau) = h(B(\Gamma),\tau)$, $A(\Gamma )$ is Koszul if and only if $A(\Gamma)^!$ is Koszul, and $A(\Gamma )^!$ is Koszul if and only if $B(\Gamma)$ is Koszul. Therefore, if either $A(\Gamma)$ or $B(\Gamma)$ is Koszul we have $h(B(\Gamma),\tau )= h(A(\Gamma), -\tau)^{-1}$.

\section{Hilbert series for $B(\Gamma)$}

\subsection{Main theorem}

\begin{theorem}
Let $\Gamma = (V, E)$ be a uniform layered graph with
$\displaystyle V = \coprod_{i=0}^N  V_i.$  Then
$$\displaystyle  h(B(\Gamma), \tau) = \displaystyle 1 + \sum_{a \in  V, |a| \ge i \ge 1} \dim \ (H^{i-2}(\Gamma_{a,i}))\tau^i.$$
\end{theorem}

\vskip 20 pt
\noindent We begin the proof with some preliminary remarks.
We will write  $V^+$  for the vector space
with basis $\{v|v \in V_+\}$. Write
$$\displaystyle {\overline v} = \sum_{w \in S(v)} w.$$
By Theorem 2.2.1, $B(\Gamma)$ has presentation
$$B(\Gamma) = T(V^+)/(I+J)$$ where
$I$ is the ideal generated by
$$\{vw|v,w \in V_+, w \notin S(v)\}$$
and $J$ is the ideal generated by $$ \{v{\overline v}| v\in V, 1 < |v|\}.$$

Let $v \in V_+$ and $n \in {\bf Z}.$
If we set $T(V^+)_{v,n}$ equal to the span of all monomials
$a_1...a_n$ with $a_1 = v$ and $a_2,...,a_n \in V_+,$ then we see that
$$\displaystyle T(V^+) = F + \oplus_{v \in V_+, 1 \leq n \leq |v| } T(V^+)_{v,n}.$$
Since  the generators of $I$ and $J$ are homogeneous with respect to this decomposition and
 $T(V^+)_{v,n} \subseteq I$ for $n > |v|$, we see that, setting
$$I_{v,n} = I \cap T(V^+)_{v,n}$$ and
$$ J_{v,n} = J \cap T(V^+)_{v,n}$$
we have
$$B(\Gamma) =
\sum_{v \in V^+, 1 \le n \le |v|} B(\Gamma)_{v,n}$$
where
$$B(\Gamma)_{v,n} = T(V^+)_{v,n}/(I_{v,n}+J_{v,n})  .$$

\vskip 15 pt

\subsection{Proof of the theorem}

Theorem 3.1.1 will follow from:
\vskip 15 pt
\begin{proposition}
$\dim(B(\Gamma)_{v,n}) = \dim(H^{n-2}(\Gamma_{v,n})).$
\end{proposition}
\vskip 15 pt
\noindent {\bf Proof:}
If $j \ge 1, u \in V_{j+2}, w \in V_j$ define
$${\overline {uw}} = \sum_{y \in S(u), w \in S(y)} y.$$

Recall that, for $0 \le j \le n-2$,
$C^{(j)}_{n-3}(\Gamma_{v,n}) =g^jC_{n-2}(\Gamma_{v,n})$ and that
$C^{(j),n-3}(\Gamma_{v,n})$ denotes the vector space of all functions
$$f:C^{(j)}_{n-3}(\Gamma_{v,n}) \rightarrow F.$$
Also write $Ch_{n-3}^{(j)}(\Gamma _{v,n})=g^jCh_{n-2}(\Gamma _{v,n})$.

Clearly
$$Ch_{n-3}(\Gamma_{v,n}) = \bigcup_{j=0}^{n-2} Ch^{(j)}_{n-3}(\Gamma_{v,n})$$
and so
$$C^{n-3}(\Gamma_{v,n}) = \bigoplus_{j=0}^{n-2} C^{(j),n-3}(\Gamma_{v,n})$$
where we extend $f \in C^{(j),n-3}(\Gamma_{v,n})$ to $Ch_{n-3}(\Gamma_{v,n})$
by setting $f(Ch^{(l)}_{n-3}(\Gamma_{v,n}))
 = 0$ for $j \neq l$.

Let $X \subseteq T(V^+)_{v,n}$ denote the span of all monomials $vb_{n-2}...b_0$ where $(b_0,...,b_{n-2}) \in Ch_{n-2}(\Gamma_{v,n}).$
Also, let
$$X_1 = span  \ \{v{\overline v}b_{n-3}...b_0|(b_0,...,b_{n-3}) \in Ch_{n-3}^{(n-2)}(\Gamma_{v,n})\},$$
$$X_j = span  \ \{vb_{n-3}...b_{n-1-j}({\overline {b_{n-1-j}b_{n-2-j}}})b_{n-2-j}...b_0|(b_0,...,b_{n-3}) \in Ch_{n-3}^ {(n-1-j)}(\Gamma_{v,n})\}$$
if $2 \leq j \leq n-2$, and
$$X_{n-1} = span  \ \{vb_{n-3}...b_0{\overline {b_0}}|(b_0,...,b_{n-3}) \in Ch_{n-3}^{(0)}(\Gamma_{v,n})\}.$$

Then we have
$$X_1,...,X_{n-1} \subseteq X,$$

$$X \cap I_{v,n} = (0),$$
$$X + I_{v,n} = T(V^+)_{v,n}$$
and
$$J_{v,n} + I_{v,n}  = I_{v,n} + \sum_{j=1}^{n-1} X_j.$$

Therefore,
$$B(\Gamma)_{v,n} = T(V^+)_{v,n}/(I_{v,n} + J_{v,n}) \cong $$
$$(X + I_{v,n})/(I_{v,n} + J_{v,n}) \cong$$
$$(X + (I_{v,n} + \sum_{j=1}^{n-1} X_j))/(I_{v,n} + \sum_{j=1}^{n-1} X_j) \cong$$
$$X/(X \cap (I_{v,n} + \sum_{j=1}^{n-1} X_j)) = X / (\sum_{j=1}^{n-1} X_j).$$

Define $$\psi: C^{n-2}(\Gamma_{v,n}) \rightarrow T(V^+)_{v,n}$$ by
 $$\psi : f \mapsto \sum_{(b_0,...,b_{n-2}) \in Ch_{n-2}(\Gamma_{v,n}) } f((b_0,...,b_{n-2}))vb_{n-2}...b_0.$$
Also, for $0 \le j \le n-2$, define
$$\psi_j: C^{(j),n-3}(\Gamma_{v,n}) \rightarrow T(V^+)_{v,n}$$
by
$$\psi_j: f \mapsto \sum_{(b_0,...,b_{n-2}) \in Ch_{n-2}(\Gamma_{v,n})} f(b_0,...,\hat b_j,...,b_{n-2})vb_{n-2}...b_0. $$

Then $\psi$ is an isomorphism of $C^{n-2}(\Gamma_{v,n})$
onto $X$, and, for $0 \le j \le n-2$,
$\psi_j$ is an isomorphism of $C^{(j),n-3}(\Gamma_{v,n})$ onto $X_{n-1-j}.$

For $f \in C^{n-3}(\Gamma_{v,n})$ write $f = \sum_{j=0}^{n-2} f_j$ with $f_j \in C^{(j),n-3}(\Gamma_{v,n}).$

Observe that
$$\psi \partial^{n-2}f = \sum_{j=0}^{n-2} (-1)^j \psi_j f_j.$$
For $0 \le j \le n-2$, define
$$C^{[j],n-3}(\Gamma_{v,n}) = \sum_{l=0}^{j} C^{(l),n-3}(\Gamma_{v,n}).$$
Note that
$$C^{[0],n-3}(\Gamma_{v,n}) =  C^{(0),n-3}(\Gamma_{v,n}) $$
and so $$C^{[0],n-3}(\Gamma_{v,n}) \cap ker \partial^{n-2}  =  (0).$$
Also note that
$$C^{[n-2],n-3}(\Gamma_{v,n}) = \sum_{l=0}^{n-2} C^{(l),n-3}(\Gamma_{v,n}) = C^{n-3}(\Gamma_{v,n}).$$
Then if $f \in C^{[j],n-3}(\Gamma_{v,n}) \cap ker \partial^{n-2}$ we have
$$0 = \sum_{l=0}^{j} (-1)^l\psi_lf_l$$
and so
$$\psi_jf_j = \sum_{l=0}^{j-1} (-1)^{l-j-1} \psi_lf_l \in X_{n-1-j} \cap(X_{n-j} + ... + X_{n-1}).$$
Furthermore, since each $\psi_l$ is surjective,
$\psi_j $ maps  $C^{[j],n-3}(\Gamma_{v,n})$  onto $X_{n-1-j} \cap(X_{n-j} + ...  + X_{n-1}).$

The kernel of the restriction of $\psi_j$
to $$C^{[j],n-3}(\Gamma_{v,n}) \cap ker \partial^{n-2}$$ is clearly
$$C^{[j-1],n-3}(\Gamma_{v,n}) \cap ker \partial^{n-2}.$$  Thus
$$\dim(X_{n-1-j} \cap (X_{n-j} + ... + X_{n-1}))  = $$
$$\dim(C^{[j],n-3}(\Gamma_{v,n}) \cap ker \partial^{n-2}) - $$
$$\dim(C^{[j-1],n-3}(\Gamma_{v,n}) \cap ker \partial^{n-2}).$$

Now we have seen that
$$ B(\Gamma)_{v,n} \cong X / (\sum_{j=1}^{n-1} X_j).$$
Since
$$\dim (\sum_{j=1}^{n-1} X_j) = \sum_{j=1}^{n-1} \dim(X_j) - \sum_{j=1}^{n-2} \dim(X_j \cap (X_{j+1} + ... + X_{n-1}))$$
we have
$$\dim(B(\Gamma)_{v,n}) = \dim(C^{n-2}(\Gamma_{v,n})) - \sum_{j=0}^{n-2} \dim (C^{(j),n-3}(\Gamma_{v,n})) + $$
$$\dim(C^{[0],n-3}(\Gamma_{v,n})\cap ker \partial^{n-2}) - \dim(C^{[n-2],n-3}(\Gamma_{v,n})\cap ker \partial^{n-2}) =$$
$$\dim(C^{n-2}(\Gamma_{v,n})) - \dim(C^{n-3}(\Gamma_{v,n})) + \dim(ker \partial^{n-2}) = $$
$$\dim(C^{n-2}(\Gamma_{v,n})) - \dim(\partial^{n-2}(C^{n-3}(\Gamma_{v,n}))).$$
Since $\partial^{n-1} = 0$ this shows that
$$\dim(B(\Gamma)_{v,n}) = \dim(H^{n-2}(\Gamma_{v,n})),$$ proving the proposition.

We may give an explicit description of the coefficients of $\tau^i$ in the Hilbert series of
$B(\Gamma)$ for $0 \le i \le 3$.  For $x \in \coprod_{i \ge 3}V_i$ write ${\mathcal E}(x,3)$
for the number of edges in the graph $\Gamma_{x,3}$ and ${\mathcal V}(x,3)$ for the number of vertices in the graph $\Gamma_{x,3}$.

\begin{corollary}
 $$h(B(\Gamma),\tau) \equiv 1 + |V_+|\tau + (\sum_{i=2}^n |E_i| - \sum_{i=2}^n |V_i|)\tau^2  + $$
$$\sum_{|x|\ge 3} ({\mathcal E}(x,3) - {\mathcal V}(x,3) + 1)\tau^3
\ mod \ (\tau^4).$$
\end{corollary}

{\bf Proof:} For any vertex $a$ and any vertex $b$ of $ \Gamma_{a,1}$ we have $a > b$ and $|a| - |b| \le 0$.  Thus the set of vertices of $\Gamma_{a,1}$ is empty and so  we have
$\dim H^{-1}(\Gamma_{a,1}) = 1$ for all $a \in V_+$.  Since $\Gamma_{a,2}$ is the graph induced by
$S(a)$ we have $\dim H^{0}(\Gamma_{a,2}) = |S(a)| - 1$ for any $a \in V_+$. Now
$\sum_{a \in V_j} |S(a)| = |E_j|$ and so
$$\sum_{a \in  V, |a| \ge 2} \dim \ (H^{0}(\Gamma_{a,2})) = \sum_{i=2}^n|E_i| - \sum_{i=2}^n|V_i|.$$
Let $x \in V_i, i \ge 3$, $u$ be a vertex of $\Gamma_{x,3}$
 of level $i-1$ and $v$ be a vertex of $\Gamma_{x,3}$
of level $i-2.$  Consider $(u)^*$ and $(v)^* \in C^0(\Gamma_{x,3}).$   We have $\partial^1(u)^* = \sum_{w \in S{(u)}} (w,u)^*$ and
$\partial^1(v^*) = -\sum_{w \in S{(x)},v \in S{(w)}} (v,w)^*.$
Now consider the $0$-chain $Y = \sum_{y \in V_1 \cup V_2} a_y(y)^*.$
The coefficient of $(v,u)^*$ in $\partial^1Y$ is $0$ unless there is an edge from $u$ to $v$ and,
if there is such an edge $e$, it is $-a_{t(e)} + a_{h(e)}.$  Thus $Y \in ker \ \partial ^1$ if and only if $a_y$ is constant on each connected component of $\Gamma_{x,3}$.
Since $\Gamma$ is uniform, $\Gamma_{x,3}$ is connected.
It follows that   $\dim H^1(\Gamma_{x,3}) = \dim C^1(\Gamma_{x,3}) - \dim C^0(\Gamma_{x,3}) +1 = {\mathcal E}(x,3) - {\mathcal V}(x,3) + 1.$

\section{$A(\Gamma)$ and numerical Koszulity}

\subsection{Hilbert series of $A(\Gamma)$}

\noindent\begin{theorem} Let $\Gamma = ({V}, {E})$ be a uniform layered graph with
$ {V} = \coprod_{i=0}^N {V}_i.$  Then:
 $$ h({A}(\Gamma),\tau)^{-1} =  1 - \sum_{a \in {V}, |a| \ge i, 0 \le s \le i-1} (-1)^s \ \dim(H^{s-1}(\Gamma_{a,i}))\tau^i.$$
\end{theorem}

{\bf Proof:}  We will use a result from \cite{RSW4}.
In Lemma 1.3 of that paper we show that if $\Gamma$ is a uniform layered graph and if, for integers $g,h$, we define
$$s_{g,h} = \sum_{v_1 > ... > v_l > *, |v_1| = g,|v_l| = h} (-1)^l,$$

then
$$ h({A}(\Gamma),\tau)^{-1} =  1 + \sum_{i\ge 1}(\sum_{g,h \in {\mathbf Z}, g \ge i \ge g-h+1} s_{g,h})\tau^i.$$

Now, for a fixed $i$, and for $a \in V$ with $|a| = g$, the chain
$$a = v_1 > ... > v_l > *$$ with $|v_l| = h$ and $|v_1| - |v_l| \le i-1$
occurs in the index set for the sum defining $s_{g,h}$ if and only if $(v_l,...,v_2) \in Ch_{l-2}(\Gamma_{a,i}).$
Thus
$$\displaystyle \sum_{g \ge i \ge g-h+1} s_{g,h} =
\displaystyle\sum_{|a| = g \ge i}\sum_{l = 1}^{i}  (-1)^{l} |Ch_{l-2}(\Gamma_{a,i})| = $$
$$\displaystyle\sum_{|a| = g \ge i}\sum_{l = 1}^{i}  (-1)^{l} \ \dim \ C^{l-2}(\Gamma_{a,i}).$$
Applying the Euler-Poincar\'e principle and setting $s = l-1$ we obtain
$$\displaystyle  -\sum_{|a| = g \ge i}\sum_{s = 0}^{i-1}  (-1)^{s} \ \dim \ H^{s-1}(\Gamma_{a,i}),$$
proving the result.

\vskip 20 pt
\subsection{Numerical Koszulity}

A version of the following theorem was announced in \cite{R}.

\begin{theorem}
 Let $\Gamma = ({V}, {E})$ be a uniform layered graph with
${V} = \coprod_{i=0}^N {V}_i.$  Assume that all minimal vertices of $\Gamma$ are contained in $V_0.$ Then
 ${A}(\Gamma)$ is numerically Koszul if and only if
$$ 0 = \sum_{a \in {V}, |a| \ge i, 0 \le s \le i-2} (-1)^s \ \dim(H^{s-1}(\Gamma_{a,i}))$$
for all $i, 3 \le i \le N.$
\end{theorem}

\noindent {\bf Proof:}  We see from Theorems 3.1.1 and 4.1.1 that $A(\Gamma)$ is numerically Koszul if and only if
$$0 = \sum_{i \ge 1}\sum_{a \in {V}, |a| \ge i, 0 \le s \le i-2} (-1)^s \ \dim(H^{s-1}(\Gamma_{a,i}))\tau^i.$$
Now the sum giving the coefficient of $\tau$ is empty. The sum giving the coefficient of $\tau^2$ is
$$\sum_{a \in V, |a| \ge 2} \ \dim ( H^{-1}(\Gamma_{a,2})).$$ Since $\Gamma_{a,2}$ is the graph induced
by the nonempty set of vertices $S(a)$,  $\dim (H^{-1}(\Gamma_{a,2})) = 0$.  Thus the coefficients of $\tau$ and $\tau^2$ are always $0$ so we have the 
result.

\vskip 20 pt
\section{Examples}
\subsection{Complete layered graphs and Boolean graphs}
We begin with two corollaries using the results of Section 1.3.

\begin{example}
The algebras $A({\bf C}[m_N,...,m_1,1])$ and $B({\bf C}[m_N,...,m_1,1])$ are numerically Koszul and
$$h(B({\bf C}[m_n,...,m_1,1]),\tau) = 1 + \sum_{k=1}^N\sum_{l=k}^N m_l(m_{l-1}-1)...(m_{l-k+1}-1)\tau^k.$$
\end{example}

\noindent {\bf Proof:} Note that if $a$ is a vertex of ${\bf C}[m_N,...,m_1,1]$ of level $j$, then
$$ {\bf C}[m_N,...,m_1,1]_{a,i} \cong {\bf C}[m_{j-1},...,m_{j-i+1}].$$
Then Proposition 1.3.2 and Theorems 3.1.1 and 4.1.1 give the result.

In fact, the algebras $A({\bf C}[m_N,...,m_1,1])$ and $B({\bf C}[m_N,...,m_1,1])$ are known to be Koszul {\cite {RSW4}} and the Hilbert series were computed in {\cite {RSW2}}.

\begin{example}
The algebras $A(\Theta_N)$ and $B(\Theta_N)$ are numerically Koszul and
$$h(B(\Theta_N),\tau) = 1 + \sum_{i=1}^N \sum_{k=i}^N {N \choose k}{{N-1} \choose {i-1}}\tau^i.$$
\end{example}
\noindent {\bf Proof:}  This follows from Theorems 3.1.1 and 4.1.1 together with the computations of Propositions 1.3.4 and 1.3.5.

In fact, the algebras $A(\Theta_N)$ and $B(\Theta_N)$ are known to be Koszul {\cite {RSW4}} and the Hilbert series were computed in {\cite {GGRSW}} and {\cite{SW}}.

\subsection{Algebras with prescribed Hilbert series}

We may use the results of this section to determine the Hilbert series of $B(\Gamma)$ for graphs $\Gamma = (V,E) $ where $V = \coprod_{i=0}^3 V_i.$

\begin{example}  Let $\Gamma = (V,E) $ be a uniform layered graph with all minimal vertices contained in $V_0$.  If $V = \coprod_{i=0}^3 V_i$,  and $V_3 = \{a\}$
and $a>b$ for any $b\in V_2$,   then $A(\Gamma)$ and $B(\Gamma)$ are Koszul and
$$h(B(\Gamma),\tau) = 1 + |V_+|\tau + (|E_2| -1)\tau^2  + (|E_2| - |V_1| - |V_2| + 1)\tau^3.$$
\end{example}

\noindent {\bf Proof:}  Since $\Gamma_{a,3} \ne \emptyset, H^{-1}(\Gamma_{a,3}) = 0.$  Then by Theorem 4.2.1, $B(\Gamma)$ is numerically Koszul if and only if
$H^0(\Gamma_{a,3}) = (0)$.  Now $ \dim \ \partial^{0}(C^{-1}(\Gamma_{a,3})) = 1$
and, by the argument of the proof of Corollary 3.2.2, $\dim \ ker(\partial^{1}(C^0(\Gamma_{a,3}))) = 1$.  Thus
$H^0(\Gamma_{a,3}) = (0)$ and $B(\Gamma)$ is numerically Koszul.  By {\cite{RSW4}}, since $V = \coprod_{i=0}^3 V_i$, the numerical Koszulity of $B(\Gamma)$ implies Koszulity.  The expression for
the Hilbert series follows from Corollary 3.2.2.

\begin{example}
Let $\Gamma$ be a graph satisfying the conditions of Example 5.2.1.  Set $r = |V_+|$ and $s = |E_2| -1.$  Then
$r \ge 3$, $r-3 \le s \le -1 + (r-1)^2/4$   and the Hilbert series of $B(\Gamma)$ is
$$h(B(\Gamma),\tau) = 1 + r\tau + s\tau^2 + (s-r+3)\tau^3.$$
Conversely, if $r,s \in \bf Z$ satisfy the above conditions, then there is a graph satisfying the conditions of Example 5.2.1 with
$r = |V_+|$ and $s = |E_2|-1.$

\end{example}

\noindent {\bf Proof:} Since $V_1,V_2,V_3 \ne \emptyset$ we have $r \ge 3$.
Since $\Gamma$ is uniform, $\Gamma_{a,3}$ must be connected and hence has at least
$|V_1| + |V_2| -1$ edges.  Of course, $\Gamma_{a,3}$ can have at most $|V_1||V_2|$ edges.  Thus
$s = |E_2|-1 \ge |V_1| + |V_2| -2 \ge |V_+| - 3 = r-3.$  Also
$s+1 \le |V_1|(r-1-|V_1|)$.  Since the maximum value of $x(r-1-x)$ is $(r-1)^2/4 $ we have the remaining inequality.
To prove the existence of such a graph, let $\Gamma' = (V',E')$ where
$V_3' = \{a\}, V_2' = \{b_1,...,b_{[(r-1)/2]}\}, V_1' = c_1,...,c_{r-1-[(r-1)/2]}\}$
and where there are edges from $a$ to every $b_i$, from $b_1$ to every $c_i$, and from every $b_i$ to $c_1$.  Then
the graph $\Gamma'_{a,3}$ is connected, so it satisfies the required conditions with $s = r-3$.  
By adding additional edges connecting vertices in $V_1$ and $V_2$ we may attain examples with $t$ edges where
$r-3 \le t \le (r-1)^2/4$ if $r$ is odd and $r-3 \le t \le r(r-2)/4$ if $r$ is even.

\begin{example}
Let $r \ge 9$. Then there is a uniform layered graph $\Gamma = (V,E)$ with $V = \coprod_{i=0}^3 V_i, V_0 = \{*\}$, and $|V_3| = 1$ such that $A(\Gamma)$ is a Koszul algebra with Hilbert series
$$(1 - r\tau + r \tau^2 - \tau^3)^{-1}.$$
\end{example}
\noindent {\bf Proof:} For example, let
$$V_3 = \{a\},$$
$$ V_2 = \{b_1,b_2,e_1,e_2,\dots , e_{r-7}\},$$
$$V_1 = \{c_1,c_2,d_1,d_2\},$$
$$E_3 = \{(a,b_1),(a,b_2)\},$$
$$ E_2 = \{(b_i,c_j),
(e_i,c_j), (e_i,d_j)\ |\ i,j = 1,2\} \cup \{(e_i,d_j)\ |\ 3 \le i \le r-7, j = 1,2\},$$
$$E_1 = \{(y,*)\ |\ y \in V_1\}$$
and apply Corollary 3.2.2.

\subsection{Algebras which are not numerically Koszul}

\begin{example} Figure 1 shows a uniform layered graph $\Gamma = (V,E)$ (due to Cassidy and Shelton) with $V = \coprod_{i=0}^4 V_i$ such that $A(\Gamma)$ is not Koszul and, in fact, is not 
numerically Koszul.  This example may be described as follows:  $$V_4 = \{a\}, V_3 = \{b_1, b_2, b_3\}, V_2 = \{c_1, c_2, c_3\}, V_1 = \{d_1, d_2, d_3\}, V_0 = \{*\};$$
$$E_4 = \{(a,b_i)|i = 1,2,3\}, E_3 = \{(b_i,c_j)|i \ne j\}, E_2 = \{(c_i,d_j)|i \ne j\}, E_1 = \{(d_i,*)|i = 1,2,3\}.$$
\newline
\begin{figure}
\centering
{\includegraphics{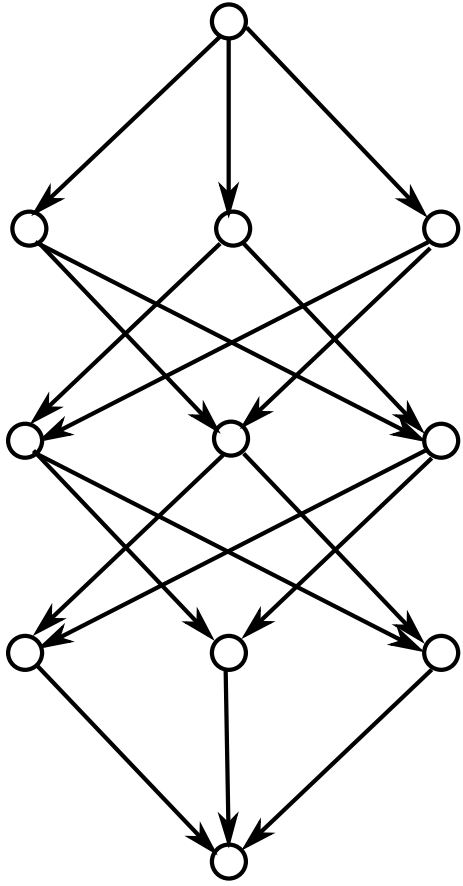}
\newline
\caption{Cassidy-Shelton graph}}
\end{figure}
\newline
\end{example}
\noindent {\bf Proof:}
We observe that
$$\dim H^{-1}(\Gamma_{a,4}) - \dim H^0(\Gamma_{a,4}) + \dim H^1(\Gamma_{a,4}) =$$
$$ \dim C^{-1}(\Gamma_{a,4}) - \dim C^0(\Gamma_{a,4}) + \dim C^1(\Gamma_{a,4}) - \dim C^2(\Gamma_{a,4}) + \dim H^2(\Gamma_{a,4}) \ge $$
$$\dim C^{-1}(\Gamma_{a,4}) - \dim C^0(\Gamma_{a,4}) + \dim C^1(\Gamma_{a,4}) - \dim C^2(\Gamma_{a,4}) =$$
$$1 - 9 + 21 - 12 = 1$$
so Theorem 4.2.1 shows that $B(\Gamma)$ is not numerically Koszul.
A similar calculation shows that the graph obtained from the Cassidy-Shelton example by deleting the edge $(b_3,c_2)$  also fails to be numerically Koszul.

\enddocument